\documentclass[12pt]{article}
\usepackage{amssymb,amsthm,bm}

\title{On meromorphic mappings admitting\\
       an Algebraic Addition Theorem}

\author{Mark B. Villarino\\
        Departamento de Matem\'atica, Universidad de Costa Rica,\\
        2060 San Jos\'e, Costa Rica}

\date{8 June 1998}

\theoremstyle{plain}
\newtheorem{thm}{Theorem}
\newtheorem{lem}{Lemma}
\newtheorem{cor}{Corollary}

\theoremstyle{definition}
\newtheorem{defn}{Definition}

\renewcommand{\theenumi}{\roman{enumi}} 
\makeatletter
\renewcommand{\p@enumii}{\theenumi}   
\@addtoreset{equation}{section}
\makeatother
 
\renewcommand{\a}{\alpha}             
\renewcommand{\b}{\beta}              
\newcommand{\f}{\varphi}              
\newcommand{\la}{\lambda}             

\newcommand{\C}{\mathbb{C}}           

\renewcommand{\aa}{{\bm{a}}}          
\newcommand{\bb}{\bm{b}}              
\newcommand{\cc}{\bm{c}}              
\newcommand{\pp}{\bm{p}}              
\newcommand{\uu}{\bm{u}}              
\newcommand{\vv}{\bm{v}}              
\newcommand{\xx}{\bm{x}}              
\newcommand{\yy}{\bm{y}}              
\newcommand{\zero}{\mathbf{0}}        

\newcommand{\A}{\mathcal{A}}          
\newcommand{\F}{\mathcal{F}}          

\newcommand{\AAT}{{\small AAT}}       

\newcommand{\Bar}{\overline}          
\newcommand{\del}{\partial}           
\newcommand{\x}{\times}               
\renewcommand{\:}{\colon}             

\newcommand{\sepword}[1]{\qquad\mbox{#1}\quad} 

\newcommand{\set}[1]{\{\,#1\,\}}      

\newcommand{\pd}[2]{\frac{\partial#1}{\partial#2}} 
\newcommand{\tpd}[2]{{\partial#1/\partial#2}} 

\newcommand{\row}[3]{{#1}_{#2},\dots,{#1}_{#3}} 

\begin{document}

\maketitle

\begin{abstract}
A proper or singular abelian mapping from $\C^n$ to $\Bar\C^n$ is 
parametrized by $n$~meromorphic functions with at most $2n$~periods. 
We develop the existence and structure theorems of the classical 
theory of an abelian mapping purely on the basis of its defining 
functional equation, the so-called algebraic addition theorem (AAT), 
with no appeal to any representation as quotients of theta functions. 
We offer two new proofs of the periodicity of a nonrational mapping 
admitting an AAT. We also prove by new arguments the existence of a 
rational group law on an associated algebraic variety, and that all
proper and singular abelian mappings do admit an~AAT.
\end{abstract}


\section{Introduction}

Weierstrass~\cite{Weierstrass1} proposed the problem of determining 
all meromorphic mappings \linebreak
$\Phi\: \C^n \to \Bar\C^n$ of complex $n$-space $\C^n$ into complex 
Osgood space $\Bar\C^n$ (the cartesian product of $n$ Riemann spheres) 
that admit an \textit{algebraic addition theorem} (\AAT). He announced 
the solution to be \textit{the set of all proper or singular abelian 
mappings}, i.e., those mappings
$$
\Phi(\uu) := \bigl( \f_1(\uu),\dots,\f_n(\uu) \bigr), 
  \qquad  \uu \in \C^n,
$$
where $\f_1(\uu),\dots,\f_n(\uu)$ are proper or singular abelian 
functions. Fifteen years later, Painlev\'e~\cite{Painleve} proved him 
right for the case $n = 2$, and in 1948 and 1954 
Severi~\cite{Severi} proved the general case.

In this paper, we obtain the fundamental properties of structure, 
periodicity and the group action of such a mapping, purely on the 
basis of its defining functional equation, the so-called \AAT, without 
ever using their explicit representation as quotients of theta 
functions, which all previous proofs have used.

We achieve this by generalizing arguments originally designed for 
elliptic functions: our means for doing so is the Weierstrass--Hurwitz 
theorem that an everywhere meromorphic function of $n$~complex
variables is rational. It does not seem to have been noticed that
this is possible.

Finally, we prove that all abelian mappings do in fact admit 
an~\AAT, by means of a general principle of algebraic dependence of 
meromorphic functions, the so-called Weierstrass--Thimm--Siegel 
theorem.

In Section~2 we explicitly state the results we will prove; subsequent 
sections contain the proofs. Section~3 is devoted to the algebraic 
dependence of the component functions and their first derivatives. In 
the following section, we construct the addition theorem variety, of 
which elliptic curves, hyperelliptic surfaces and Picard varieties are 
particular instances. Section~5 contains two new proofs that an 
abelian mapping is either rational or periodic, as well as a vivid 
interpretation of what periodicity is. In Section~6, we give a proof 
the addition theorem can be expressed rationally in terms of $n + 1$ 
meromorphic functions on the addition theorem variety. In the next 
section, we use this rational addition theorem to define an abelian 
group law on this variety, which is the starting point 
of the modern theory of abelian varieties. 

Finally, in Section~8, we  give a direct proof that all \emph{abelian}
mappings admit an  algebraic addition theorem; therefore, all the
previous results are applicable to such mappings. We emphasize that
we assume only that our mapping is meromorphic and admits an \AAT, and
then obtain as a \emph{theorem} the classical definition, that it is
parameterized by rational or periodic functions with at most
$2n$~periods.

\section{Statement of results}

We begin with the following fundamental definition. Let 
$\f_1(\uu),\dots,\f_n(\uu)$ be $n$ ana\-lytically independent 
meromorphic functions on $\C^n$. They define a meromorphic mapping 
$\Phi\: \C^n \to \Bar\C^n$ where
$$
\Phi : \uu \mapsto (\f_1(\uu),\dots,\f_n(\uu)).
$$

\begin{defn}
\label{df:AAT}
We say that $\Phi$ \textbf{admits an algebraic addition theorem
(\AAT)} if and only if there exist $n$ polynomials
$$
G_k \equiv G_k(\la; \row{x}{1}{n}, \row{y}{1}{n}),  \qquad 
 k = 1,\dots,n,
$$
in $(2n + 1)$ variables $\la$ and $\row{x}{1}{n}$ and $\row{y}{1}{n}$ 
with complex coefficients, such that the equations 
\begin{equation}
G_k\{ \f_k(\uu + \vv); \f_1(\uu),\dots,\f_n(\uu),
      \f_1(\vv),\dots,\f_n(\vv) \} = 0,  \qquad  k = 1,\dots,n
\label{eq:aat-Gk}
\end{equation}
hold, and none of the denominators vanish identically in~$\uu$ 
and~$\vv$. 
\end{defn}

In this paper, we prove the following fundamental properties of the 
mapping $\Phi$ \textit{purely on the basis of the \AAT}, without 
making use of the results of Weierstrass, Painlev\'e and Severi on its 
explicit analytic form.

\begin{thm}
\label{th:alg-dep}
{\rm(a)} Any $(n + 1)$ of the $(n^2 + n)$ functions
$$
\f_1(\uu),\dots,\f_n(\uu),\sepword{and} \pd{\f_k}{u_p}(\uu), \quad
 k,p = 1,\dots,n,
$$
are \textbf{algebraically dependent}.

{\rm(b)}
Any function $\f_k(\uu)$ is \textbf{algebraically dependent} on its 
$n$ partial derivatives of first order, $\tpd{\f_k}{u_p}(\uu)$, 
$p = 1,\dots,n$.

{\rm(c)}
Any one of the $n^2$ \textbf{first-order partial derivatives} 
$\tpd{\f_k}{u_p}(\uu)$ $(k,p = 1,\dots,n)$ is algebraically dependent 
on the $n$ original functions $\f_1(\uu),\dots,\f_n(\uu)$.
\end{thm}

\begin{thm}
\label{th:aat-variety}
The $n^2$ partial derivatives $\tpd{\f_k}{u_p}(\uu)$, for 
$k,p = 1,\dots,n$, generate a \textbf{simple algebraic extension} 
$\F$ of the field of rational functions of $\f_1(\uu),\dots,\f_n(\uu)$
with complex coefficients. The minimal polynomial
$V(\theta;\row{x}{1}{n})$  of this extension is satisfied by
$$
x_1 := \f_1(\uu),\ \dots, \ x_n := \f_n(\uu), \quad
\theta := \sum_{k,p=1}^n \a_{kp} \pd{\f_k}{u_p}(\uu),
$$
where, for suitable $\a_{kp} \in \C$, the function $\theta$ is a 
\textbf{primitive element} of~$\F$. Thus
$$
V\{\theta(\uu); \f_1(\uu),\dots,\f_n(\uu)\} = 0.
$$
\end{thm}

\begin{defn}
\label{df:aat-variety}
The hypersurface in $\Bar\C^{n+1}$ defined by
\begin{equation}
V(\theta;\row{x}{1}{n}) = 0
\label{eq:aat-variety}
\end{equation}
is called the \textbf{addition theorem variety}~$\A$.
\end{defn}

\begin{thm}
\label{th:periodicity}
Any generator $\row{\f}{1}{n}$ of~$\F$ that is not rational is 
\textbf{periodic}.
\end{thm}

\begin{thm}[{\rm Rational form of the \AAT}]
\label{th:rational-aat}
Let $\uu \in \C^n$ and $\vv \in \C^n$ be two independent variables 
and let $\theta,\row{\f}{1}{n}$ generate the addition theorem 
variety~$\A$. Then there exist $(n + 1)$ \textbf{rational} functions
$$
R_k(x_0,\row{x}{1}{n}; y_0,\row{y}{1}{n}),  \qquad  k = 0,\dots,n
$$
of $(2n + 2)$ variables $\row{x}{0}{n},\row{y}{0}{n}$ with constant 
coefficients, such that the equations
$$
\f_k(\uu + \vv) = R_k\{ \theta(\uu), \f_1(\uu),\dots,\f_n(\uu);
                        \theta(\vv), \f_1(\vv),\dots,\f_n(\vv) \}
$$
for $k = 0,1,\dots,n$ hold, and none of the denominators vanish 
identically in $\uu$ and $\vv$.
\end{thm}

\begin{thm}
\label{th:group-law}
The rational addition theorem (Theorem~\ref{th:rational-aat}) defines 
an \textbf{abelian group law} on the addition theorem variety~$\A$.
\end{thm}

\begin{thm}
\label{th:abelian-aat}
Every proper or singular abelian mapping admits an algebraic addition 
theorem.
\end{thm}

\section{Algebraic dependence of first order derivatives}

Theorem~\ref{th:alg-dep}a is the general result, while 
parts~\ref{th:alg-dep}b and~\ref{th:alg-dep}c are immediate 
corollaries. We prove this theorem by means of an explicit elimination 
process which amounts to a finite algorithm that uses only rational 
operations. Painlev\'e~\cite{Painleve} obtains related results by a 
completely different procedure that involves the reversion of infinite 
series, something which our algorithm avoids completely, and he does 
not obtain the general principle stated in Theorem~\ref{th:alg-dep}a.

Theorems~\ref{th:alg-dep}b and~\ref{th:alg-dep}c are well-known 
properties of the abelian functions, but our proof shows that 
\textit{any} meromorphic mapping that has an \AAT\ enjoys these two 
properties, and that the \AAT\ is the primordial reason for that. 
Thus, all \textit{rational} mappings have these properties, since they 
admit an~\AAT.

\begin{proof}[Proof of Theorem~\ref{th:alg-dep}]
The \AAT~(\ref{eq:aat-Gk}) for the mapping~$\Phi$ can be rewritten as
\begin{equation}
G_k[ \f_k(\uu + \vv); \Phi(\uu), \Phi(\vv) ] = 0,  \qquad 
 k = 1,\dots,n,
\label{eq:aat-GkPhi}
\end{equation}
where $G_1(\la_1;\xx,\yy) \not\equiv 0$, \dots, 
$G_n(\la_n;\xx,\yy) \not\equiv 0$ are polynomials in the $2n + 1$ 
variables $\xx = (\row{x}{1}{n})$, $\yy = (\row{y}{1}{n})$ and~$\la_k$. 
Differentiating 
\begin{equation}
G_1[ \f_1(\uu + \vv); \Phi(\uu), \Phi(\vv) ] = 0
\label{eq:aat-G1Phi}
\end{equation}
with respect to~$u_1$ and then $v_1$, we get
\begin{eqnarray}
\pd{G_1}{\la_1} \pd{\la_1}{u_1}
 + \sum_{k=1}^n \pd{G_1}{x_k} \pd{x_k}{u_1} & = & 0,
\label{eq:dG1-du1}
\\
\pd{G_1}{\la_1} \pd{\la_1}{v_1}
 + \sum_{k=1}^n \pd{G_1}{y_k} \pd{y_k}{v_1} & = & 0.
\label{eq:dG1-dv1}
\end{eqnarray}
Now
\begin{equation}
\pd{\la_1}{u_1} = \pd{\f_1(\uu+\vv)}{u_1}
 = \pd{\f_1(\uu+\vv)}{(u_1+v_1)}
 = \pd{\f_1(\uu+\vv)}{v_1} = \pd{\la_1}{v_1}.
\label{eq:equal-derivs}
\end{equation}
Therefore, subtracting (\ref{eq:dG1-dv1}) from (\ref{eq:dG1-du1}) and 
using~(\ref{eq:equal-derivs}), we arrive at
\begin{equation}
\sum_{k=1}^n \biggl( \pd{G_1}{x_k} \pd{x_k}{u_1}
  - \pd{G_1}{y_k} \pd{y_k}{v_1} \biggr) = 0,
\label{eq:delG1-u1v1}
\end{equation}
whose left hand side is a polynomial in the variables 
$\row{x}{1}{n}$, $\row{y}{1}{n}$, $\tpd{x_1}{u_1}$, \dots, 
$\tpd{x_n}{u_1}$, $\tpd{y_1}{v_1}$, \dots, $\tpd{y_n}{v_1}$ 
and~$\la_1$. If we take the greatest common divisor 
of~(\ref{eq:delG1-u1v1}) and the \AAT~(\ref{eq:aat-G1Phi}), we 
obtain an equation of the form
$$
g_{11}\biggl(\la_1; \row{x}{1}{n}, \pd{x_1}{u_1},\dots,\pd{x_n}{u_1},
       \row{y}{1}{n}, \pd{y_1}{v_1},\dots,\pd{y_n}{v_1} \biggr) = 0,
$$
and the remainder, set equal to zero, is the \emph{eliminant}:
$$
H_{11}\biggl(\row{x}{1}{n}, \pd{x_1}{u_1},\dots,\pd{x_n}{u_1},
        \row{y}{1}{n}, \pd{y_1}{v_1},\dots,\pd{y_n}{v_1} \biggr) = 0.
$$

If we now differentiate (\ref{eq:delG1-u1v1}) with respect to~$u_p$ 
and then $v_p$ for $p = 2,\dots,n$, we get the greatest common 
divisors
$$
g_{1p}\biggl(\la_1; \row{x}{1}{n}, \pd{x_1}{u_p},\dots,\pd{x_n}{u_p},
       \row{y}{1}{n}, \pd{y_1}{v_p},\dots,\pd{y_n}{v_p} \biggr) = 0,
$$
and the eliminants
$$
H_{1p}\biggl(\row{x}{1}{n}, \pd{x_1}{u_p},\dots,\pd{x_n}{u_p},
        \row{y}{1}{n}, \pd{y_1}{v_p},\dots,\pd{y_n}{v_p} \biggr) = 0.
$$
We now apply the same process to the \AAT s $G_2 = 0$, $G_3 = 0$, 
\dots, $G_n = 0$. Then we get a total of $n^2$ greatest common divisors
$$
g_{kp}\biggl(\la_k; \row{x}{1}{n}, \pd{x_1}{u_p},\dots,\pd{x_n}{u_p},
       \row{y}{1}{n}, \pd{y_1}{v_p},\dots,\pd{y_n}{v_p} \biggr) = 0
$$
and $n^2$ eliminant equations:
$$
H_{kp}\biggl(\row{x}{1}{n}, \pd{x_1}{u_p},\dots,\pd{x_n}{u_p},
        \row{y}{1}{n}, \pd{y_1}{v_p},\dots,\pd{y_n}{v_p} \biggr) = 0,
$$
for $k,p = 1,\dots,n$. We observe that each polynomial $g_{kp}$ and 
$H_{kp}$ is symmetric in $(x_1,y_1)$, \dots, $(x_n,y_n)$, 
$(\tpd{x_1}{u_p}, \tpd{y_1}{v_p})$, \dots, 
$(\tpd{x_n}{u_p}, \tpd{y_n}{v_p})$ as a consequence of the 
\AAT~(\ref{eq:aat-GkPhi}) and the symmetry of $(\uu + \vv)$ in~$\uu$ 
and~$\vv$.

If we suitably \emph{fix} the variable~$\vv$, the $n^2$ eliminant 
equations become
\begin{equation}
h_{kp}\biggl(\row{x}{1}{n}, \pd{x_1}{u_p},\dots,\pd{x_n}{u_p} \biggr)
 = 0,  \qquad  k,p = 1,\dots,n,
\label{eq:elims-hkp}
\end{equation}
where $h_{kp}$ is a polynomial in the $2n$ variables 
$\row{x}{1}{n}$, $\tpd{x_1}{u_p},\dots,\tpd{x_n}{u_p}$. Therefore the 
$n^2$ equations~(\ref{eq:elims-hkp}) relate the $n + n^2$ variables 
$\row{x}{1}{n}$, $\tpd{x_k}{u_p}$ ($k,p = 1,\dots,n$), and we can 
eliminate any $(n^2 - 1)$ of them, leaving the remaining $(n + 1)$ of 
them as algebraically dependent. This proves Theorem~\ref{th:alg-dep}.
\end{proof}

\section{The addition theorem variety}

Theorem~\ref{th:aat-variety} reveals the origin of elliptic curves, 
elliptic hypersurfaces, \dots, abelian varieties. They are all 
manifestations of the \textit{addition theorem variety}~$\A$, 
generated by the $n^2$ partial; derivatives $\tpd{\f_k}{u_p}(\uu)$ 
$(k,p = 1,\dots,n)$ over the field of rational functions of 
$\f_1(\uu),\dots,\f_n(\uu)$. Painlev\'e~\cite{Painleve} also arrives 
at~$\A$, although in not quite so explicit a form, and from a 
slightly different point of view. He transforms the $n^2$ equations 
for the first order partial derivatives into a system of total 
(algebraic) differential equations
\begin{eqnarray}
du_1 
&=& p_{11}(\row{x}{1}{n})\,dx_1 +\cdots+ p_{1n}(\row{x}{1}{n})\,dx_n,
\nonumber \\
\vdots && \qquad \vdots
\nonumber \\
du_n 
&=& p_{n1}(\row{x}{1}{n})\,dx_1 +\cdots+ p_{nn}(\row{x}{1}{n})\,dx_n,
\label{eq:painleve-sys}
\end{eqnarray}
where $x_1 = \f_1(\uu)$, \dots, $x_n = \f_n(\uu)$, and
$$
p_{ij} := \pd{u_i}{x_j} = \frac{1}{J} \pd{J}{z_{ij}},  \qquad  
 (i,j = 1,\dots,n),
$$
where $z_{ij} := \tpd{\f_i}{u_j}(\uu)$ and $J$ is the determinant of 
the Jacobian matrix $[z_{ij}]$, $i,j = 1,\dots,n$. (We have altered 
Painlev\'e's notation to conform with ours, and have considered the 
case of general~$n$, instead of $n = 2$ as he does.) Now he states:

\begin{quote}
``\dots\ as is well known, one may rationally express [the algebraic 
irrationalities] $p_{ij}(\row{x}{1}{n})$ in terms of $\row{x}{1}{n}$ 
and a unique irrationality $\theta(\row{x}{1}{n})$ defined by means 
of an algebraic relation 
$$
S(\theta;\row{x}{1}{n}) = 0,
$$
which is such that conversely, $\theta$ is rationally expressible in 
$\row{x}{1}{n}$, $p_{ij}$ $(i,j = 1,\dots,n)$. Since $p$ or 
$\tpd{u_1}{x_1}$ is rationally deducible from $\tpd{u_i}{x_j}$
$(i,j = 1,\dots,n)$, as also are [all the other] $p_{ij}$, the 
function $\theta(\uu)$ is uniform and meromorphic at the same time as 
$x_1(\uu),\dots,x_n(\uu)$ \cite[\S 7]{Painleve}.''
\end{quote}

He never explicitly constructs $\theta$ as 
$\sum_{k,p=1}^n \a_{kp}z_{kp}$ with $\a_{kp} \in \C$, which is crucial 
for our further development of the theory. His entire approach is 
based on the system~(\ref{eq:painleve-sys}) on the addition theorem 
variety~$\A$, which he calls ``the algebraic surface parametrized by 
the functions $\row{x}{1}{n}$''.

\begin{proof}[Proof of Theorem~\ref{th:aat-variety}]
This is an immediate consequence of elementary field theory, the 
primitive element theorem~\cite[p.~243]{Lang} and
Theorem~\ref{th:alg-dep}c, the corollary that every first-order
partial derivative $\tpd{\f_k}{u_p}(\uu)$ is algebraically dependent 
on $\f_1(\uu),\dots,\f_n(\uu)$.
\end{proof}

\section{Periodicity}

We present \textit{two proofs} that the mapping $\Phi$ is periodic. 
Our proofs are based on the Weierstrass--Hurwitz theorem that a 
meromorphic function in~$\Bar\C^n$ with no essential singularities is 
a \textit{rational} function of all its variables, and vice versa; 
and on a weak form of Picard's theorem on the behaviour of an 
analytic function of one variable in the vicinity of an isolated 
essential singularity. Both proofs show that the ``cause'' or 
``explanation'' of the existence of a period $\pp \neq 0$ of the 
meromorphic mapping~$\Phi$ is \textit{the simultaneous occurrence of 
two antithetical properties} of~$\Phi$:
\begin{itemize}
\item the wild chaotic dispersion of values of~$\Phi$ in the 
neighbourhood of the essential singularity at infinity (an 
``irresistible force'');
\item the rigid unyielding restriction on the values of $\Phi(\uu)$ 
imposed by the polynomial form of the \AAT\ (an ``immovable object'').
\end{itemize}
The mathematical resolution of this ancient philosophical conundrum is 
this: \textit{the solution set $\uu \in \C^n$ of the equation 
$\Phi(\uu) = \xx$, $\xx \in \Bar\C^n$ fixed, is a discrete 
$k$-dimensional real lattice}, $0 \leq k \leq 2n$. This beautiful 
interpretation of periodicity is due, in principle, to Weierstrass for 
the case of functions of one variable.

Painlev\'e~\cite{Painleve} also proves that $\Phi$ is periodic. But 
his proof is totally different. He starts from the 
system~(\ref{eq:painleve-sys}) of total (algebraic) differential 
equations and shows that $\row{u}{1}{n}$ are line integrals on~$\A$ 
of the first, second or third kinds. Then the topology of~$\A$ shows 
that its fundamental group is generated by $k \leq 2n$ linearly 
independent cycles. Our proofs, on the other hand, avoid the topology 
of~$\A$ and the theory of line integrals on~$\A$, and work directly 
with the polynomials of the \AAT\ and the singularities of~$\Phi$ at 
infinity.

\begin{defn}
Let $\Phi$ be the meromorphic mapping
$$
\Phi \: \C^n \to \Bar\C^n
      : \uu \mapsto \bigl( \f_1(\uu),\dots,\f_n(\uu) \bigr)
$$
defined by the $n$ meromorphic functions $\row{\f}{1}{n}$. We say 
that $\Phi$ is \textbf{periodic} if and only if there exists a 
nonzero $\pp \in \C^n$ such that
$$
\Phi(\uu + \pp) = \Phi(\uu)
$$
for all $\uu \in \C^n$ where $\Phi$ is defined. The vector~$\pp$ is 
called a \textbf{period} of~$\Phi$.
\end{defn}

We shall present two proofs of Theorem~\ref{th:periodicity} that are 
based on the existence of an essential singularity plane at infinity 
of at least one of the component functions 
$\f_1(\uu),\dots,\f_n(\uu)$. The first proof is direct, while the 
second uses an important sufficient condition for periodicity.

\subsection{Picard's ``tiny'' theorem}

Picard's ``great'' theorem affirms that if a meromorphic function 
$\f(n)$ of one complex variable~$u$ has an essential singularity at 
infinity, then for any choice of complex number~$c$ (with at most two 
exceptions), the solution set of the equation $\f(u) = c$ is an 
infinite discrete set on the complex $u$-plane. In his lectures, 
Weierstrass proved the following ``tiny'' version at least twenty 
years before Picard:

\begin{thm}
Suppose that the function $\f(u)$ has an isolated singularity at the 
point $u = a$. Let $c$ be an arbitrary complex constant and let 
$|w - c| < h$ be an arbitrarily small neighbourhood of the point~$c$. 
Then there exists a point~$c'$ in this neighbourhood such that the 
equation $\f(u) = c'$ has \emph{infinitely many roots} which 
accumulate onto the point~$a$.
\end{thm}

The elegant elementary proof can be found in Hancock~\cite{Hancock}, 
Osgood~\cite{Osgood} and Phrag\-m\'en~\cite{Phragmen}.

\subsection{The solution set $\Phi(\uu) = \aa$}

The analytic hypersurface $S_\aa$ defined by the vector equation
$$
\Phi(\uu) = \aa,  \qquad  \aa \in \C^n \mbox{ fixed},
$$
is a complicated set of points. However, we are interested in a very 
elementary property of~$S_\aa$.

\begin{thm}
\label{th:infty-solns}
Suppose that $\Phi(\uu)$ is \emph{not} a rational mapping, i.e., that 
at least one of the coordinate functions $\f_k(\uu)$ ($k = 1,\dots,n$) 
is transcendental. Suppose further that $\Phi$ is holomorphic at the 
point $\uu = \aa$ and that $\Phi(\uu) = \aa$. Then the set $S_\aa$ 
contains infinitely many distinct points $\uu_n$ that accumulate on at 
least one plane at infinity.
\end{thm}

\begin{proof}
We assume that the Jacobian determinant is nonzero in a neighbourhood 
of $\uu = \aa := (\row{a}{1}{n})$. Then the set of equations becomes
$$
\f_1(\uu) = a_1, \ \f_2(\uu) = a_2, \ \dots\ \f_n(\uu) = a_n.
$$
Suppose that $\f_n(\uu)$ is transcendental; then it is not a rational 
function of $\row{u}{1}{n}$. Now the Weierstrass--Hurwitz 
theorem~\cite{Hurwitz}
states that $\f_n(\uu)$ is rational in~$\uu$ if and only if it has no 
essential singularities in~$\Bar\C^n$. But $\f_n(\uu)$, being 
transcendental, does have essential singularities; however, it is 
meromorphic in all of~$\C^n$. Therefore it has at least one plane at 
infinity as an essential singularity surface. Suppose that plane is 
$u_n = \infty$.

We now revert the equations
$$
\f_1(\uu) = a_1, \ \dots\ \f_{n-1}(\uu) = a_{n-1}
$$
to get
\begin{equation}
u_1 = F_1(u_n), \ u_2 = F_2(u_n), \ \dots\ u_{n-1} = F_{n-1}(u_n)
\label{eq:uj-Fun}
\end{equation}
and substitute these expansions into $\f_n(\uu) = a_n$, to get an 
equation of the form
$$
\f_n\{F_1(u_n),\dots,F_{n-1}(u_n),u_n\} = a_n,
$$
or $F_n(u_n) = a_n$ (say), which is a meromorphic function in 
finite~$\C$ with an essential singularity at infinity. Therefore, by 
Picard's theorem, there is an infinity of values 
$$
u_n = a_0, v_1, \row{v}{2}{k},\dots
$$
which are distinct roots of $F_n(u_n) = a_n$ and accumulate on 
$u_n = \infty$. Any one of these, say $u_n = v'$, when substituted 
into~(\ref{eq:uj-Fun}) will give values $u_j = v_j$ 
($j = 1,\dots,n-1$) such that
$$
\f_1(\row{v}{1}{n-1},v') = a_1, \ \f_{n-1}(\row{v}{1}{n-1},v') = a_{n-1}, 
 \ \dots\ \f_n(\row{v}{1}{n-1},v') = a_n.
$$
Thus, as $u_n = v'$ runs through the infinitely many roots of 
$F_n(u_n) = a_n$, the points $(\row{v}{1}{n-1},v')$ run through 
infinitely many values that satisfy $\Phi(\uu) = \aa$ and accumulate 
on $u_n = \infty$.
\end{proof}

\subsection{First proof of periodicity}

\begin{proof}[Proof of Theorem~\ref{th:periodicity}]
The \AAT\ for a particular $\f := \f_k$ has the form
\begin{equation}
G[\f(\uu + \vv); \Phi(\uu), \Phi(\vv)] = 0.
\label{eq:aat-phi}
\end{equation}
Let $\cc_2 := \Phi(\vv)$ for some $\vv \in \C^n$. Then, by 
Theorem~\ref{th:infty-solns},
$$
V := \set{\vv \in \C^n : \Phi(\vv) = \cc_2}
$$
contains infinitely many distinct elements $\vv_1,\vv_2,\dots$ that 
accumulate on at least one of the planes at infinity. Let $m$ be the 
degree of~$G$ in the first variable $\f(\uu + \vv)$ and choose 
$m + 1$ values of~$\vv_k$. We assume that we choose $\uu \in \C^n$ 
such that 
$$
\uu,\ \uu + \vv_1,\ \uu + \vv_2,\ \dots,\ \uu + \vv_{m+1}
$$
are nonsingular points of~$\Phi$; this we can do since there are only
a finite number of the~$\vv_k$. Then, if
\begin{equation}
\cc_1 := \Phi(\uu),
\label{eq:c1-phiu}
\end{equation}
the equation~(\ref{eq:aat-phi}) becomes
$$
G[\f(\uu + \vv_k); \cc_1, \cc_2] = 0,
$$
which is an algebraic equation of degree~$m$ whose roots are the 
$m + 1$ numbers
$$
\f(\uu + \vv_1),\ \f(\uu + \vv_2),\ \dots,\ \f(\uu + \vv_{m+1}),
$$
and so at least two are equal:
\begin{equation}
\f(\uu + \vv_k) = \f(\uu + \vv_l).
\label{eq:equal-roots}
\end{equation}
This holds for any~$\uu$ in a neighbourhood of the original one. 
Moreover, since $\f$ is holomorphic in a neighbourhood of~$\uu$ and 
$\Phi$ is holomorphic in the translated points $\uu + \vv_k$, then 
$\Phi$ is holomorphic in a neighbourhood of each $\uu + \vv_k$ and we 
can take each such neighbourhood to be the translate of a sufficiently 
small neighbourhood of~$\uu$. For each~$\uu$, we get a 
different~$\cc_1$ in~(\ref{eq:c1-phiu}), and so a possibly different 
pair $\vv_k$ and~$\vv_l$ in~(\ref{eq:equal-roots}). But there are only 
$(m + 1)^2$ pairs $(\vv_k,\vv_l)$ available, so there is an entire 
neighbourhood of values of~$\uu$ in which (\ref{eq:equal-roots})~holds 
for the \emph{same} pair $(\vv_k,\vv_l)$. This means that
$$
\f(\uu + \vv_k) \equiv \f(\uu + \vv_l)
$$
identically; or equivalently, $\f(\uu + \vv_k - \vv_l)\equiv \f(\uu)$,
and so $\vv_k - \vv_l$ is a period of~$\f$.
\end{proof}

\subsection{Second proof of periodicity}

\subsubsection{The algebraic differential equations}

The periodicity of at least one of the functions 
$\f_1(\uu),\dots,\f_n(\uu)$ depends on the existence of algebraic 
differential equations of the first order. In 
Theorem~\ref{th:alg-dep}c we proved that each one of the $n^2$ 
first-order partial derivatives $\tpd{\f_k}{u_p}(\uu)$ 
$(k,p = 1,\dots,n)$ is algebraically dependent on these $n$ original 
functions; that is, there exist $n^2$ polynomials 
$P_{kp}(z_0;\row{z}{1}{n})$ in the $(n + 1)$ variables 
$z_0,\row{z}{1}{n}$ with complex coefficients such that the following 
relations hold identically:
\begin{equation}
P_{kp}\biggl(\pd{\f_k}{u_p}(\uu);\f_1(\uu),\dots,\f_n(\uu)\biggr)
 \equiv 0.
\label{eq:pddep-phi}
\end{equation}

\begin{thm}
\label{th:rat-dep-derivs}
Suppose that the analytically independent functions $\row{\f}{1}{n}$ 
admit an \AAT. Then the partial derivatives of orders two and beyond 
are uniquely determined as rational functions of the $n^2$ first-order 
partial derivatives $\tpd{\f_k}{u_p}$ $(k,p = 1,\dots,n)$.
\end{thm}

\begin{proof}
By Theorem~\ref{th:alg-dep}c, there exist $n^2$ polynomials 
$P_{kp}(z_0;\row{z}{1}{n})$ satisfying~(\ref{eq:pddep-phi}). 
Let us abbreviate $P_{kp,i} = \tpd{P_{kp}}{z_i}$ for 
$i = 0,1,\dots,n$. Differentiating the equations $P_{kp} = 0$ with 
respect to~$u_q$, we obtain
\begin{equation}
\frac{\del^2\f_k}{\del u_p\del u_q}
 = - \frac{1}{P_{kp,0}} \biggl( P_{kp,1} \pd{\f_1}{u_q}
        +\cdots+ P_{kp,n} \pd{\f_n}{u_q} \biggr).
\label{eq:rat-dep-derivs}
\end{equation}
It is obvious how to continue to derivatives of higher order.
Therfore, all partial derivatives of the $P_{kp}$, i.e., all the
$P_{kp,i}$, all the $P_{kp,ij} = \del^2 P_{kp}/\del z_i\del z_j$,
and so on, are polynomials in $\row{\f}{1}{n}$ and $\tpd{\f_i}{u_q}$ 
$(i,q = 1,\dots,n)$ whose coefficients are complex constants.
\end{proof}

\subsubsection{A sufficient condition for periodicity}

The following theorem generalizes to $n$ complex variables a 
sufficient condition for periodicity first stated and proved (in 
lectures) by Weierstrass. As we shall see, it is an essential element 
in the theory of functions that admit an~\AAT.

Recall that the \emph{derivative} $\Phi'(\uu)$ of the mapping~$\Phi$ 
is the $n \x n$ Jacobian matrix whose $(k,p)$-entry is the partial 
derivative $\tpd{\f_k}{u_p}(\uu)$.

\begin{thm}
\label{th:one-period}
Suppose that the meromorphic mapping~$\Phi$ admits an \AAT. If there 
exist two distinct points $\aa \in \C^n$ and $\bb \in \C^n$ such that 
\begin{equation}
\Phi(\aa) = \Phi(\bb)  \sepword{and}\quad  \Phi'(\aa) = \Phi'(\bb),
\label{eq:phiata-phiatb}
\end{equation}
and if each $P_{kp,0} \neq 0$, then the mapping $\Phi$ is periodic 
with period $\pp := \aa - \bb$.
\end{thm}

\begin{proof}
By assumption, $\f_1(\uu),\dots,\f_n(\uu)$ are holomorphic in a 
neighbourhood of both $\aa$ and~$\bb$. By Taylor's theorem,
\begin{eqnarray}
\f_k(\uu + \aa) 
&=& \f_k(\aa) + \sum_{p=1}^n \pd{\f_k}{u_p}(\aa)\,u_p + \frac{1}{2} 
 \sum_{p,q=1}^n \frac{\del^2\f_k}{\del u_p\del u_q}(\aa) \,u_p\,u_q 
 + \cdots
\\
\f_k(\uu + \bb) 
&=& \f_k(\bb) + \sum_{p=1}^n \pd{\f_k}{u_p}(\bb)\,u_p + \frac{1}{2}
 \sum_{p,q=1}^n \frac{\del^2\f_k}{\del u_p\del u_q}(\bb) \,u_p\,u_q 
 + \cdots
\label{eq:taylor-phi}
\end{eqnarray}
Using the notation
$$
\a_k := \f_k(\aa) = \f_k(\bb),  \quad 
\a_{k,p} := \pd{\f_k}{u_p}(\aa) = \pd{\f_k}{u_p}(\bb),
$$
and obvious extensions to higher derivatives, the 
equations~(\ref{eq:taylor-phi}) become
\begin{eqnarray}
\f_k(\uu + \aa) 
&=& \a_k + \sum_{p=1}^n \a_{k,p}\,u_p
         + \frac{1}{2} \sum_{p,q=1}^n \a_{k,pq} \,u_p \,u_q  + \cdots
\\
\f_k(\uu + \bb) 
&=& \a_k + \sum_{p=1}^n \a_{k,p}\,u_p
         + \frac{1}{2} \sum_{p,q=1}^n \b_{k,pq} \,u_p \,u_q  + \cdots
\label{eq:taylor-phi-abrev}
\end{eqnarray}
But by Theorem~\ref{th:rat-dep-derivs} and the 
assumption~(\ref{eq:phiata-phiatb}), all the higher partial 
derivatives of $\row{\f}{1}{n}$ are uniquely determined by the values 
of the $\a_k$ and the~$\a_{k,p}$; which means [and 
here we use $P_{kp,0} \neq 0$, when plugging 
in~(\ref{eq:rat-dep-derivs})]:
$$
\a_{k,pq} = \b_{k,pq}, \ \a_{k,pqr} = \b_{k,pqr}, \ \dots  
 \sepword{for all} k,p,q,r,\dots
$$
and therefore the two expansions~(\ref{eq:taylor-phi-abrev}) are 
identical. This means that the equation
$$
\Phi(\uu + \aa) \equiv \Phi(\uu + \bb)
$$
holds identically in~$\uu$, which implies that 
$\Phi(\uu + \aa - \bb) \equiv \Phi(\uu)$, and that therefore 
$\pp := \aa - \bb$ is a period vector.
\end{proof}

\subsubsection{Completion of the proof}

\begin{proof}[Second proof of Theorem~\ref{th:periodicity}]
We assume that $\f_1(\uu),\dots,\f_n(\uu)$ are $n$ analytically 
independent meromorphic functions that admit an \AAT. Moreover, we 
assume that $\f_n(\uu)$ is \emph{not} a rational function of~$\uu$. We 
know the set $S_\aa = \set{\uu \in \C^n : \Phi(\uu) = \aa}$ contains 
an infinite number of points $\vv_k$ such that 
$v_{kn} \to \infty$ as $k \to \infty$.

We now apply a Dirichlet pigeon-hole principle argument to the
$n^2$-tuples of values $\tpd{\f_k}{u_p}(\uu)$ $(k,p = 1,\dots,n)$ as
$\uu$ runs over the points $\vv_k \in S_\aa$. Each such function
value is a root of the algebraic equations
$$
P_{kp}\biggl(\pd{\f_k}{u_p}(\uu);\row{a}{1}{n}\biggr) = 0,  \quad
 (k,p = 1,\dots,n),
$$
of respective degrees $m_{kp}$ in the first variables 
$\tpd{\f_k}{u_p}$. Therefore, the maximum number of distinct 
$n^2$-tuples is the product $\prod_{k,p=1}^n m_{kp} < \infty$. 
Therefore, as $\uu$ runs over the infinite set of pairs $\{\vv_k\}$, 
there is an infinite subsequence $\bb_1,\row{\bb}{2}{r},\dots$ 
in~$S_\aa$ such that
$$
\pd{\f_k}{u_p}(\bb_1) = \pd{\f_k}{u_p}(\bb_2) = \cdots 
 \sepword{for each}  k,p = 1,\dots,n.
$$
These equations more than fulfil the 
sufficient conditions of Theorem~\ref{th:one-period} for~$\f_n$ to be 
periodic with the period vector $\pp = \bb_1 - \bb_2$. Thus 
we have proved that $\f_n$ is a periodic function.
\end{proof}

\section{The rational addition theorem}

The only published proof of the rational form of the AAT 
(Theorem~\ref{th:rational-aat}), is due to 
Siegel~\cite[pp.~94--96]{Siegel3}, and is only for the case in which 
$\f_1(\uu),\dots,\f_n(\uu)$ are $n$ independent abelian functions. We 
shall adapt it to the case of \emph{any} mapping~$\Phi$ with an~\AAT.
Let $\f$ be any one of the~$\f_k$. Then Siegel's proof consists 
of the following steps:
\begin{enumerate}

\item  Let $\F$ be the field of abelian functions generated by
       $\f_1(\uu),\dots,\f_n(\uu)$. Then $\F$ is a \emph{simple 
	   algebraic extension} of $\C(\f_1(\uu),\dots,\f_n(\uu))$ with 
	   primitive element $\f_0(\uu)$.
\label{it:alg-extn}

\item  For each fixed $\vv = \bb$, the function $\f(\uu + \bb)$ is an 
       element of~$\F$, and is therefore a \emph{rational} function
       of $\f_0(\uu),\f_1(\uu),\dots,\f_n(\uu)$.
\label{it:phiub-rat}

\item  For each fixed $\vv = \bb$, the function $\f(\uu + \bb)$ 
       belongs to~$\F$, and for each fixed $\uu = \aa$, the function 
	   $\f(\aa + \vv)$ belongs to~$\F$. Therefore $\f(\uu + \vv)$ lies 
	   in $\F \otimes \F$, i.e., it is rational \emph{jointly} in 
	   $\f_0(\uu),\f_1(\uu),\dots,\f_n(\uu)$ and
	   $\f_0(\vv),\f_1(\vv),\dots,\f_n(\vv)$.
\label{it:phiuplusv-rat}
\end{enumerate}

There is a large literature on the question of proving that functions 
that are rational in each variable separately are also rational in
those variables jointly; see, e.g.,
\cite{Andrea,Hurwitz,Osgood,Palais} among others. Siegel uses the
argument expounded in the book of Bochner and Martin
\cite[pp.~199--203]{BochnerM} (wherein further sources are cited),
although he does not say so explicitly.

Siegel proves (\ref{it:alg-extn}) by a detailed analysis of the 
period matrix and the ``Thetasatz'' of Weierstrass and 
Riemann~\cite{Siegel3}, which affirms that every abelian function 
belonging to a given period matrix can be represented as a quotient of 
two (generalized) theta functions, to prove that 
$\f_0(\uu),\f_1(\uu),\dots,\f_n(\uu)$ satisfy a polynomial equation 
$$
P\bigl[ \f_0(\uu),\f_1(\uu),\dots,\f_n(\uu) \bigr] = 0,
$$
where $P(x_0,\row{x}{1}{n}) \not\equiv 0$ is a polynomial with 
constant coefficients such that
\begin{enumerate}
\item[I.] its degree $q$ in $x_0$ does not depend on~$\f_0$, but is 
       \emph{uniquely determined} by $\row{\f}{1}{n}$;
\item[II.] its degree $r$ in $(\row{x}{1}{n})$ does not depend
       on~$\uu$.
\end{enumerate}
Therefore any $\f_0$ whose degree in~$P$ is maximal will be a 
primitive element for~$\F$.

In our case, we define $\F$ \emph{directly} as 
$$
\F := \C(\row{x}{1}{n})(\theta),
$$
and try to adapt Siegel's reasoning to our situation. Thus, 
step~(\ref{it:alg-extn}) for us is trivial, since it's true by 
definition.

Step~(\ref{it:phiub-rat}) is trivial for Siegel, since if $\f(\uu)$ is 
an abelian function so is $\f(\uu + \bb)$, because it belongs to the 
same period matrix and is also meromorphic. However, 
step~(\ref{it:phiub-rat}) is not trivial for us. We prove it in two 
stages:
\begin{enumerate}
\item[(a)] For each fixed~$\vv = \bb$, the function $\f(\uu + \bb)$ 
           lies in a \emph{finite algebraic extension} of 
			$\C(\f_1(\uu),\dots,\f_n(\uu))$ whose degree does not
			depend on~$\bb$.
\item[(b)] The degree of this extension is \emph{one}, and therefore 
           $\f(\uu + \bb)$ lies in~$\F$.
\end{enumerate}
The assertion~(a) is an immediate consequence of the \AAT, while 
for~(b) it suffices to prove:
\begin{enumerate}
\item[(b$'$)] The function $\f(\uu + \bb)$ is \emph{uniquely 
            determined} by $(\theta(\uu),\f_1(\uu),\dots,\f_n(\uu))$; 
			that is, given any $(n+1)$-tuple 
			$(\theta;\row{x}{1}{n}) \in\A$, any $\uu\in\C^n$ such that
			$$
			\theta(\uu) = \theta; \ \f_1(\uu) = x_1, \ \dots, \ 
			\f_n(\uu) = x_n
			$$
			will give \emph{the same numerical value} of 
			$\f(\uu + \bb)$.
\end{enumerate}

Our proof of~(b$'$) uses our explicit construction of the primitive 
element~$\theta$, the separability of the minimal polynomial 
$V(\theta;\row{x}{1}{n})$ of the extension $\F : \C(\row{x}{1}{n})$, 
of which $\theta$ and all its conjugates are the roots, and the 
sufficient condition, Theorem~\ref{th:one-period}, for the periodicity 
of~$\Phi$.

Step~(\ref{it:phiuplusv-rat}) is not trivial for Siegel nor for us.
But the fundamental lemma which he proves, on the basis of 
properties~(I) and~(II) of his 
polynomial $P(\row{x}{0}{n})$, are trivial for us and our polynomial 
$V(\theta;\row{x}{1}{n})$. To make the argument work, Siegel proves the 
following lemma.

\begin{lem}
For each $k = 0,1,\dots,h$, there exists an algebraic relation between
$$
\theta_k := \f(\uu + \vv)\, \theta^{h-k}
$$
and $\row{x}{1}{n}$, whose total degree does not exceed a bound~$r$ 
independent of~$\vv$.
\end{lem}

\begin{proof}
(We have altered Siegel's notation of to conform with ours.) Siegel 
proves this lemma by appealing to the ``Thetasatz'' of Weierstrass and 
Riemann~\cite{Siegel3}, which affirms that every abelian function 
belonging to a given period matrix can be represented as a quotient of 
two (generalized) theta functions.

But this claim is trivial from our point of view, since the degree of 
the \AAT\ for $\f(\uu + \vv)$ does not depend on~$\vv$, nor does the 
degree of the field polynomial $V(\theta;\row{x}{1}{n})$ defining the 
algebraic extension~$\F$. The degree of~$V$ in \emph{any} of its 
$(n + 1)$ variables does not depend on~$\vv$. But the polynomial 
over $\C(\row{x}{1}{n})$ whose roots are the $\theta_k$ is obtainable 
from the coefficients of the polynomial giving the \AAT\ for $\f(\uu)$ 
and the polynomial $V(\theta;\row{x}{1}{n})$ by means of a finite 
number of rational operations, as given in the theory of elimination. 
Therefore, the total degree will be bounded by an integer~$r$ which 
does not depend on~$\vv$.
\end{proof}

Once we get past that step, Siegel's proof carries over word for word, 
and we are done.

\begin{proof}[Proof of Theorem~\ref{th:rational-aat}]
The proof consists in verifying the following four lemmas.
\end{proof}

\begin{lem}
\label{lm:one}
For each fixed~$\vv = \bb$, the function $\f_k(\uu + \bb)$ is an 
element of a finite algebraic extension of~$\F$, for $k = 1,\dots,n$.
\end{lem}

\begin{proof}
The \AAT\ for $\f_k(\uu + \vv)$ is of the form~(\ref{eq:aat-Gk}), 
where $G_k(\la;\xx,\yy)$ is a polynomial with constant coefficients, 
not all zero, in the $(2n + 1)$ variables 
$\la,\row{x}{1}{n},\row{y}{1}{n}$. Taking $\vv = \bb$ gives
$$
G_k\bigl[ \f_k(\uu + \bb); \Phi(\uu),\Phi(\bb) \bigr] = 0,
$$
and this is a polynomial whose coefficients are rational functions of 
$\f_1(\uu),\dots,\f_n(\uu)$ and whose roots are the values of 
$\f_k(\uu + \bb)$. That is, they lie in a finite extension field of 
$\C(\row{x}{1}{n})$ and thus of $\C(\theta; \row{x}{1}{n})$ also. 
\end{proof}

\begin{lem}
\label{lm:two}
Each value of the primitive element $\theta$ of~$\F$ uniquely 
determines the $n^2$ values of the first-order partial derivatives 
$\tpd{\f_k}{u_p}(\uu)$, for $k,p = 1,\dots,n$.
\end{lem}

\begin{proof}
By definition,
$$
\theta = \sum_{k,p=1}^n \a_{kp}\, \pd{\f_k}{u_p}(\uu),
$$
for suitable $\a_{kp} \in \C$. As each of the partial derivatives 
$\tpd{\f_k}{u_p}(\uu)$ independently runs over all its respective 
conjugate values, the function $\theta$ takes on all of its 
$h$~distinct values as the distinct roots of the minimal polynomial
\begin{equation}
V(\theta; \row{x}{1}{n}) = 0.
\label{eq:min-poly}
\end{equation}
Thus, each choice of $n^2$~values of $\tpd{\f_k}{u_p}(\uu)$ gives a 
unique $\theta$ and each $\theta$ gives a unique set of $n^2$
partial derivatives. 
\end{proof}

\begin{lem}
\label{lm:three}
Each $(n+1)$-tuple $(\theta;\row{x}{1}{n}) \in \A$, where $\A$ is 
the addition theorem variety, \emph{uniquely determines} the value 
of $\f_k(\uu + \bb)$.
\end{lem}

\begin{proof}
Fix a particular point $(\theta; \row{x}{1}{n}) \in \A$. The the $n$ 
numbers $\row{x}{1}{n}$ determine a solution set $U$ of the equations
\begin{equation}
\f_1(\uu) = x_1, \ \dots, \ \f_n(\uu) = x_n,
\label{eq:phiu-x}
\end{equation}
that is, $U = \set{\uu \in \C^n : \Phi(\uu) = \xx}$. Let $|U|$ denote 
the cardinality of~$U$, and write 
$\Theta := \set{\theta(\uu) : \uu \in U}$. 

We consider three possibilities. Firstly, if $|U|$ is less than the 
degree~$h$ in~$\theta$ of the polynomial~$V$ in~(\ref{eq:min-poly}), 
then $|\Theta| < h$, so there are fewer than~$h$ values~$\theta$, 
which is impossible.

Secondly, if $|U| = h$, then one value of~$\theta$ corresponds to one 
value of~$\uu$ and vice versa, so that $\uu$ is uniquely determined, 
and therefore $\f_k(\uu + \bb)$ is uniquely determined, since 
$\f_k$ is single-valued.

Lastly, if $|U| > h$, then $|\Theta| > h$, so that more than $h$ 
values of~$\theta$ satisfy the field equation~(\ref{eq:min-poly}), of 
degree~$h$ in the $theta$~variable. Therefore, $\theta(\uu)$ takes on 
the same value at two distinct members $\uu = \aa_1,\aa_2 \in U$. Now 
$\Phi'(\aa_1) = \Phi'(\aa_2)$ on account of Lemma~\ref{lm:two}, and 
$\Phi(\aa_1) = \Phi(\aa_2)$ also by assumption~(\ref{eq:phiu-x}). By 
Theorem~\ref{th:one-period}, $\Phi(\uu)$ is periodic with period 
$\pp = \aa_1 - \aa_2$. Hence $\f_k(\aa_1 + \bb) = \f_k(\aa_2 + \bb)$ 
for any~$\bb$. The conclusion is that 
$\f_k(\uu_1 + \bb) = \f_k(\uu_2 + \bb)$ whenever $\uu_1,\uu_2 \in U$ 
satisfy $\theta(\uu_1) = \theta(\uu_2)$, so that in this case also,
$\f_k(\uu + \bb)$ is uniquely determined.
\end{proof}

\begin{lem}
\label{lm:four}
For each fixed~$\vv = \bb$, the function $\f_k(\uu + \bb)$ is an 
element of the field~$\F$.
\end{lem}

\begin{proof}
By Lemma~\ref{lm:one}, $\f_k(\uu + \bb)$ lies in an extension field 
of~$\F$, and by Lemma~\ref{lm:three}, $\f_k(\uu + \bb)$ is a 
\emph{single-valued} algebraic function of $(\theta;\row{x}{1}{n})$, 
which means that the degree of the polynomial defining 
$\f_k(\uu + \bb)$ in an extension field of~$\F$ is of \emph{degree 
one}, that is to say, $\f_k(\uu + \bb)$ lies in~$\F$ itself.
\end{proof}

\section{The group law}

Our point of view leads to a very explicit form of the group law on 
the addition theorem variety $\A$ \emph{in terms of the functions and 
their first derivatives}. This explicit form has not been cited in the 
literature, much less proved, for the case of $n > 1$ variables. Of 
course, it is a well-known theorem of Liouville in the case $n = 1$ of 
elliptic functions. Our development shows that our form of the group 
law is the ``natural'' one that springs organically from the explicit 
polynomials of the~\AAT.

The converse is also true. That is, we have proved that if $\Phi$ 
admits an~\AAT, then it determines an algebraic variety with an 
abelian group law. Picard~\cite{Picard} proved that if an algebraic 
variety has a transitive group law defined on it, then the variety 
can be parametrized by abelian functions, and it is well known that 
they have an~\AAT. Moreover, if the group is not transitive, then the 
variety can be parametrized by singular abelian functions, and they 
in turn can be proved to admit an \AAT\ (see 
Section~\ref{sc:sing-abel}). This last result is the crowning 
achievement of half a century of effort, starting with the memoir of 
Painlev\'e~\cite{Painleve} and culminating in the work of 
Severi~\cite{Severi}.

\begin{lem}
If $\f(\uu)$ belongs to~$\F$, so also does $\f(-\uu)$.
\end{lem}

\begin{proof}
Apply the \AAT~(\ref{eq:aat-Gk}) of $\Phi(\uu + \vv)$ to 
$\Phi(\uu - \uu)$:
$$
G_k[\f_k(\zero); \Phi(\uu),\Phi(-\uu)] = 0, \qquad k = 1,\dots,n,
$$
which can be reexpressed as
\begin{equation}
D_k[\f_1(\uu),\dots,\f_n(\uu); \f_1(-\uu),\dots,\f_n(-\uu)] = 0, 
 \qquad  k = 1,\dots,n,
\label{eq:dk-poly}
\end{equation}
for some polynomial $D_k(\row{x}{1}{n};\row{y}{1}{n}) \not\equiv 0$.
(Note that if $D_k$ were to vanish identically, we could use instead 
$\uu = \uu_1 + \bb$, $\vv = -\uu_1 - \bb$ for a suitable value 
of~$\bb$.) If we eliminate, say, $\f_2(-\uu),\dots,\f_n(-\uu)$ 
from~(\ref{eq:dk-poly}), we obtain
$$
E_1[\f_1(-\uu); \f_1(\uu),\dots,\f_n(\uu)] = 0, 
$$
which means that $\f_1(-\uu)$ is algebraically dependent on 
$\f_1(\uu),\dots,\f_n(\uu)$. The same holds for 
$\f_2(-\uu),\dots,\f_n(-\uu)$. Thus $\Phi(-\uu)$ admits an \AAT\ 
wherein the coefficients of the powers of $\f_k(-\uu + \vv)$ are 
rational functions of $\Phi(\uu)$ and~$\Phi(\vv)$. The reasoning we 
used for the rational addition theorem now shows that 
$\f(-\uu + \vv)$ lies in~$\F$ for each fixed~$\vv$. In particular, 
$\f(-\uu)$ lies in~$\F$.
\end{proof}

\begin{cor}[{\rm Subtraction theorem}]
$\f_k(\uu - \vv)$ lies in $\F \otimes \F$.
\qed
\end{cor}

\begin{proof}[Proof of Theorem~\ref{th:group-law}]
We define the operation of \textbf{addition} on the addition theorem 
variety~$\A$ as follows. If
$$
P_1 := \bigl( \theta(\uu); \f_1(\uu),\dots,\f_n(\uu) \bigr)
\quad\mbox{and}\quad
P_2 := \bigl( \theta(\vv); \f_1(\vv),\dots,\f_n(\vv) \bigr),
$$
then
$$
P_1 \oplus P_2 := \bigl(\theta(\uu + \vv);
                          \f_1(\uu + \vv),\dots,\f_n(\uu + \vv)\bigr),
$$
and we define \emph{subtraction} by
$$
P_1\ominus P_2 := \bigl(\theta(\uu - \vv);
                          \f_1(\uu - \vv),\dots,\f_n(\uu - \vv)\bigr).
$$
This evidently defines an abelian group on~$\A$ with identity 
$\bigl( \theta(\zero); \f_1(\zero),\dots,\f_n(\zero) \bigr)$. In fact, 
for each fixed~$\vv$, the rational addition theorem, 
Theorem~\ref{th:rational-aat}, yields a birational regular mapping 
of~$\A$ onto itself, and therefore we obtain a \emph{group} of 
birational regular mappings of the addition theorem variety onto 
itself.
\end{proof}

\section{Proper and singular abelian functions}
\label{sc:sing-abel}

We can reformulate Theorem~\ref{th:periodicity} as follows:
\begin{quote}
\emph{A meromorphic mapping $\Phi$ that admits an algebraic addition 
theorem is either 
(1)~a \textbf{rational} mapping, or 
(2)~a \textbf{periodic} mapping.}
\end{quote}

\begin{defn}
If $\Phi\: \C^n \to \Bar\C^n$ is periodic, and if each component 
function $\f_k$ of $\Phi = (\row{\f}{1}{n})$ has $\pi_k$ periods, then
\begin{enumerate}
\item[(a)] If $\pi_k = 2n$ for all~$k$, then $\Phi$ is a 
           \textbf{proper abelian mapping}.
\item[(b)] If $\pi_k < 2n$ for at least one~$k$, then $\Phi$ is a
           \textbf{singular abelian mapping}.
\end{enumerate}
The component $\f_k$ is called a \emph{proper} or \emph{singular} 
abelian \emph{function}, according as $\pi_k = 2n$ or $\pi_k < 2n$.
\end{defn}

There is no standard terminology about the case $\pi_k < 2n$. Here is 
a list of terms used by various authors:
\begin{center}
\begin{tabular}{l@{\qquad}l}
Weierstrass & Degenerate abelian function;\\
Painlev\'e  & Degenerate abelian function;\\
Cousin      & Semirational abelian function~\cite{Cousin};\\
Severi      & Quasi-abelian function;\\
Siegel      & Singular abelian function.
\end{tabular}
\end{center}
The term ``degenerate'' is a good one since such functions arise when 
the faces of the period parallelotope of a proper abelian function are 
translated to infinity (this last fact is quite difficult to prove). 
However, the term \emph{degenerate} is used nowadays for meromorphic 
functions whose period group group is not a lattice, something quite 
different from its meaning for Weierstrass. We have chosen Siegel's 
term ``singular'' as expressing the same property.

Note that $\Phi$ is proper if and only if all its component functions 
are proper. $\Phi$ is singular if and only if \emph{two} conditions 
hold:
\begin{enumerate}
\item  $\pi_k < 2n$ for some~$k$;
\item  $\Phi$ admits an \AAT.
\end{enumerate}

Our proof of Theorem~\ref{th:abelian-aat} is based on a famous result 
that has fascinated mathematicians for more than a century, namely the 
Weierstrass--Thimm--Siegel theorem on algebraic dependence.

\begin{thm}[{\rm Weierstrass-Thimm-Siegel}]
\label{th:WThS}
Let $\row{\f}{1}{k}$ be meromorphic functions on a compact complex 
space~$X$. Then they are \emph{algebraically dependent} if
\begin{enumerate}
\item  they are analytically dependent; or if
\label{it:first}
\item  $k = n + 1$ and $X$ has dimension~$n$; or if
\label{it:second}
\item  $X$ is irreducible and $\row{\f}{1}{k}$ are analytically 
       dependent on a nonempty open subset of~$X$.
\label{it:third}
\end{enumerate}
\end{thm}

The main result is (\ref{it:first}), while (\ref{it:second}) and 
(\ref{it:third}) are corollaries. Weierstrass announced 
(\ref{it:second}) for the case that $\row{\f}{1}{n+1}$ are abelian 
functions belonging to the period parallelotope 
$X$~\cite{Weierstrass2}, although he never published a complete proof. 
The first complete proof of~(\ref{it:first}) was given by 
Thimm~\cite{ThimmThesis} in his K\"onigsberg thesis in~1939, and then 
Siegel gave two proofs, in 1948~\cite{SiegelA} and 
1955~\cite{SiegelM}. Thimm's historical review~\cite{ThimmW} is an 
excellent source on the origins and proofs of this theorem. Suffice 
it to say that nowadays the proof is so elementary that it has 
appeared in standard textbooks~\cite{Shafarevich,Siegel3}.

This wonderful theorem makes the proof of Theorem~\ref{th:abelian-aat} 
for the case of proper abelian mappings almost trivial, and has 
nothing with the explicit analytical form of the functions~$\f_k$.

The case that $\Phi$ is a \emph{singular} abelian mapping presents 
significant difficulties not found in the case of proper mappings. 
The problem is that in the latter case, the period parallelotope is a 
\emph{compact} complex manifold, while there is no obvious compact 
counterpart for the singular mappings. Even in the case $n = 1$, the 
function $e^{\sin u}$ is simply periodic and indeed entire on~$\C$, 
but it does not admit an~\AAT; its domain of definition is a 
noncompact period strip instead of a period parallelogram.

If $n = 1$, the singular elliptic functions turn out to be
\begin{enumerate}
\item  rational functions of $u$; and
\item  rational functions of $e^{au}$, for $a \in \C$.
\end{enumerate}
If $n = 2$, the singular abelian functions turn out to be
\begin{enumerate}
\item  rational functions of $\f_1 = u_1$ and $\f_2 = u_2$;
\item  rational functions of $\f_1 = u_1$ and $\f_2 = e^{u_2}$;
\item  rational functions of $\f_1 = e^{u_1}$ and $\f_2 = e^{u_2}$;
\item  rational functions of $\f_1 = \wp(u_1)$ and 
       $\f_2 = u_2 - \epsilon\, \zeta(u_1)$;
\label{it:fourth}
\item  rational functions of $\f_1 = \wp(u_1)$ and 
       $\f_2 = e^{u_2}\,\sigma(u_1 - a)/\sigma(u_1)$;
\end{enumerate}
where $\epsilon = 0$ or~$1$, $\zeta(u_1)$ is the Weierstrass 
$\zeta$-function, $a$ is arbitrary and $\sigma(u_1)$ is the Weierstrass 
$\sigma$-function. That these are singular abelian functions is easy 
to see; Painlev\'e's great achievement was to prove that \emph{there
are no others}.

As an example, consider case~(\ref{it:fourth}) with $\epsilon = 1$. 
Here
$$
\f_1(u_1,u_2) = \wp(u_1),  \qquad  \f_2(u_1,u_2) = u_2 - \zeta(u_1).
$$
It is a standard result that $\zeta(u_1)$ is quasiperiodic, that is, if 
$2\omega_1$ and $2\omega_2$ are the periods of $\wp(u_1)$, then
$$
\zeta(u_1 + 2\omega_1) = \zeta(u_1) + 2\eta_1,  \qquad  
\zeta(u_1 + 2\omega_2) = \zeta(u_1) + 2\eta_2,
$$
where $\eta_1 := \zeta(\omega_1)$, $\eta_2 := \zeta(\omega_2)$. 
Therefore
$$
\f_2(u_1 + 2\omega_1, u_2 + 2\eta_1) 
 = u_2 + 2\eta_1 - \zeta(u_1 + 2\omega_1)
 = u_2 + 2\eta_1 - [\zeta(u_1) + 2\eta_1] = \f_2(u_1,u_2),
$$
and similarly $\f_2(u_1 + 2\omega_2, u_2 + 2\eta_2) = \f_2(u_1,u_2)$. 
This means that
$$
\pp_1 := (2\omega_1, 2\eta_1)  \sepword{and}\quad
\pp_2 := (2\omega_2, 2\eta_2) 
$$
are the two periods of the singular mapping $\Phi = (\f_1,\f_2)$.

We are left with the problem of characterizing the singular abelian 
mappings intrinsically, so as to be able to prove 
Theorem~\ref{th:abelian-aat}. The cases $n = 1$ and $n = 2$ show that 
such functions can be regarded as \emph{rational functions of a mobile 
point $(\row{x}{1}{p},\row{y}{1}{q})$ on the cartesian product}
$$
\A = \A_p \x \Bar\C^q,
$$
where $n = p + q$, $\A_p$ is a Picard variety of dimension~$p$ 
defined by the equation
$$
V(\theta; \row{x}{1}{p}) = 0,
$$
and $(\row{y}{1}{q})$ are the coordinates of a point of~$\C^q$. We 
shall \emph{define} a singular abelian function of $n$~complex 
variables to be such a rational function. Note that this is the
definition that Severi uses in his great memoir~\cite{Severi}.

\begin{proof}[Proof of Theorem~\ref{th:abelian-aat}]
Case~1: suppose that $\Phi$ is a \emph{proper} abelian mapping. 
Then $\f_k(\uu + \vv)$, for each $k = 1,\dots,n$, is a meromorphic 
function defined on the $2n$-dimensional compact complex manifold 
$X \x X$ where $X$ is the period parallelotope for $\Phi$. By 
Theorem~\ref{th:WThS},
\begin{equation}
\f_k(\uu+\vv);\ \f_1(\uu),\dots,\f_n(\uu),\ \f_1(\vv),\dots,\f_n(\vv),
\label{eq:alg-dep}
\end{equation}
being $(2n + 1)$ meromorphic functions on the $2n$-dimensional compact 
complex manifold $X \x X$, are \emph{algebraically dependent}; that is 
to say, the mapping $\Phi$ admits an~\AAT.

\vspace{6pt}

Case~2: $\Phi$ is a \emph{singular} abelian mapping.
Now $\f_k(\uu + \vv)$, for each $k = 1,\dots,n$, is a meromorphic 
function defined on the $2n$-dimensional compact complex manifold 
$(\A_p \x \Bar\C^q) \x (\A_p \x \Bar\C^q)$, where $\A_p$ is a Picard 
variety and $p + q = n$. Thus the functions~(\ref{eq:alg-dep}) are 
again algebraically dependent, by the Weierstrass--Thimm--Siegel 
theorem.
\end{proof}

\subsection*{Acknowledgments}
I are most grateful to the International Centre for Theoretical 
Physics at Trieste for its hospitality, stimulating atmosphere and 
library facilities; the idea for this paper germinated during a recent 
visit there. I also thank R.~Balasubramanian for several conversations 
and insightful criticisms, and especially Joseph~C. V\'arilly for 
many discussions and unrelenting encouragement. Support from the 
Vicerrector\'{\i}a de Investigaci\'on of the University of Costa Rica 
is acknowledged.


\end{document}